# Some integrals involving the Stieltjes constants

Donal F. Connon

dconnon@btopenworld.com

12 February 2009

**Abstract**

Some new integrals involving the Stieltjes constants are developed in this paper.

**1. Introduction**

Using the binomial theorem we obtain

$$\int_0^\infty e^{-ux}(1-e^{-x})^n x^{p-1} dx = \int_0^\infty e^{-ux} \sum_{j=0}^n \binom{n}{j}(-1)^j e^{-jx} x^{p-1} dx$$

$$= \sum_{j=0}^n \binom{n}{j}(-1)^j \int_0^\infty e^{-(u+j)x} x^{p-1} dx$$

and employing the definition of the gamma function for $\text{Re}(p) > 0$

$$\Gamma(p) = \int_0^\infty e^{-x} x^{p-1} dx$$

we see for $\text{Re}(u) > 0$ that

$$\frac{\Gamma(p)}{u^p} = \int_0^\infty e^{-ux} x^{p-1} dx$$

Hence for $\text{Re}(u) > 0$ and $\text{Re}(p) > 0$ we have

(1.1) $$\int_0^\infty e^{-ux}(1-e^{-x})^n x^{p-1} dx = \Gamma(p) \sum_{j=0}^n \binom{n}{j} \frac{(-1)^j}{(u+j)^p}$$

With $n = 0$ in (1.1) we have $\int_0^\infty e^{-ux} x^{p-1} dx = \frac{\Gamma(p)}{u^p}$ and therefore (1.1) is also valid for $n = 0$. The integral (1.1) appears in Gradshteyn and Ryzhik [16, p.357, Eq. 3.432.1] where it is attributed to Bierens de Haan.

As $p \to 0$ the right-hand side of (1.1) is indeterminate but we may use L'Hôpital's rule when considering the ratio

$$\frac{\sum_{j=0}^{n}\binom{n}{j}\frac{(-1)^j}{(u+j)^p}}{1/\Gamma(p)}$$

(since $\sum_{j=0}^{n}\binom{n}{j}(-1)^j = \delta_{n,0}$ we require that $n \neq 0$ in order to ensure that the numerator is equal to zero as $p \to 0$).

As regards the numerator, it is easily seen that

$$\frac{d}{dp}\sum_{j=0}^{n}\binom{n}{j}\frac{(-1)^j}{(u+j)^p} = -\sum_{j=0}^{n}\binom{n}{j}\frac{(-1)^j \log(u+j)}{(u+j)^p}$$

We have using Euler's reflection formula for the gamma function

$$\frac{d}{dp}[1/\Gamma(p)] = \frac{1}{\pi}\frac{d}{dp}[\Gamma(1-p)\sin \pi p]$$

$$= \Gamma(1-p)\cos \pi p - \frac{1}{\pi}\Gamma'(1-p)\sin \pi p$$

and therefore we see that where $r$ is a positive integer

$$\lim_{p \to 1-r}\frac{d}{dp}[1/\Gamma(p)] = (-1)^{r+1}\Gamma(r)$$

This gives us

(1.2) $$\int_0^\infty \frac{e^{-ux}(1-e^{-x})^n}{x^r}dx = \frac{(-1)^r}{\Gamma(r)}\sum_{j=0}^{n}\binom{n}{j}(-1)^j(u+j)^{r-1}\log(u+j)$$

The above integral was originally given by Anglesio [3] in 1997. [Could one adopt a Eulerian approach at this stage and assume that $r$ was a continuous variable with the representation $(-1)^r = \cos \pi r$ ?]

We now let $u = 1$ and $r = 1$ to obtain

(1.3) $$\int_0^\infty \frac{e^{-x}(1-e^{-x})^n}{x}\,dx = -\sum_{j=0}^{n}\binom{n}{j}(-1)^j \log(1+j)$$

Making the substitution $y = e^{-x}$, equation (1.3) becomes for $n \geq 1$

(1.4) $$\int_0^1 \frac{(1-y)^n}{\log y}\,dy = \sum_{j=0}^{n}\binom{n}{j}(-1)^j \log(1+j)$$

and this is in agreement with G&R [16, p.541, Eqn. 4.267.1] as corrected by Boros and Moll (see the Errata to G&R [16] dated 26 April 2005).

Making the same substitution $y = e^{-x}$ in (1.2) we see that for $n \geq 1$

(1.5) $$\int_0^1 \frac{y^{u-1}(1-y)^n}{\log^r y}\,dy = \frac{1}{\Gamma(r)}\sum_{j=0}^{n}\binom{n}{j}(-1)^j (u+j)^{r-1}\log(u+j)$$

and we have the summation (where we are required to start at $n = 1$)

$$\sum_{n=1}^{\infty}\frac{1}{n+1}\int_0^1 \frac{y^{u-1}(1-y)^n}{\log^r y}\,dy = \frac{1}{\Gamma(r)}\sum_{n=1}^{\infty}\frac{1}{n+1}\sum_{j=0}^{n}\binom{n}{j}(-1)^j (u+j)^{r-1}\log(u+j)$$

Using the logarithmic expansion

$$\sum_{n=1}^{\infty}\frac{(1-y)^n}{n+1} = -\frac{1}{1-y}[\log y + 1 - y]$$

we therefore obtain

(1.6) $$\int_0^1 \left[\frac{1}{1-y} + \frac{1}{\log y}\right]\frac{y^{u-1}}{\log^{r-1} y}\,dy = -\frac{1}{\Gamma(r)}\sum_{n=1}^{\infty}\frac{1}{n+1}\sum_{j=0}^{n}\binom{n}{j}(-1)^j (u+j)^{r-1}\log(u+j)$$

Starting the summation at $n = 0$ gives us

(1.7)

$$\int_0^1 \left[\frac{1}{1-y} + \frac{1}{\log y}\right]\frac{y^{u-1}}{\log^{r-1} y}\,dy = \frac{1}{\Gamma(r)}u^{r-1}\log u - \frac{1}{\Gamma(r)}\sum_{n=0}^{\infty}\frac{1}{n+1}\sum_{j=0}^{n}\binom{n}{j}(-1)^j (u+j)^{r-1}\log(u+j)$$

and reference to the Hasse identity [18]

$$(s-1)\varsigma(s) = \sum_{n=0}^{\infty} \frac{1}{n+1} \sum_{j=0}^{n} \binom{n}{j} \frac{(-1)^j}{(u+j)^{s-1}}$$

shows that the right-hand side of (1.7) is related to $\varsigma'(2-r,u)$.

With $r=1$ we get

(1.8) $$\int_0^1 \left[\frac{1}{1-y} + \frac{1}{\log y}\right] y^{u-1} dy = \log u - \sum_{n=0}^{\infty} \frac{1}{n+1} \sum_{j=0}^{n} \binom{n}{j} (-1)^j \log(u+j)$$

and it is known from [17] and [12] that the digamma function $\psi(u)$ may be represented by

(1.9) $$\psi(u) = \sum_{n=0}^{\infty} \frac{1}{n+1} \sum_{j=0}^{n} \binom{n}{j} (-1)^j \log(u+j)$$

Equation (1.8) then becomes

(1.10) $$\int_0^1 \left[\frac{1}{1-y} + \frac{1}{\log y}\right] y^{u-1} dy = \log u - \psi(u)$$

This integral is due to Binet [22, p.175] and with $u=1$ in (1.10) we obtain the well-known integral for Euler's constant [6, p.178]

(1.11) $$\gamma = \int_0^1 \left[\frac{1}{1-y} + \frac{1}{\log y}\right] dy$$

The identity (1.9) may also be obtained in the following way by employing a method due to Dirichlet [25]. We see that

$$\frac{\Gamma(u+h) - \Gamma(u)}{h} = \frac{\Gamma(u+h)}{\Gamma(1+h)} \Gamma(h) - \frac{\Gamma(u+h)}{\Gamma(1+h)} \frac{\Gamma(u)\Gamma(h)}{\Gamma(u+h)}$$

$$= \frac{\Gamma(u+h)}{\Gamma(1+h)} [\Gamma(h) - B(u,h)]$$

$$= \frac{\Gamma(u+h)}{\Gamma(1+h)} \int_0^1 \left[\frac{1}{\log^{1-h}(1/y)} - \frac{y^{u-1}}{(1-y)^{1-h}}\right] dy$$

and in the limit as $h \to 0$ we have

$$\Gamma'(u) = -\Gamma(u)\int_0^1 \left[\frac{1}{\log y} + \frac{y^{u-1}}{1-y}\right]dy$$

Therefore we have

(1.12) $$\psi(u) = -\int_0^1 \left[\frac{1}{\log y} + \frac{y^{u-1}}{1-y}\right]dy$$

We see that

$$\frac{1}{u} = \int_0^1 y^{u-1}\,dy$$

and integration of this results in the well-known Frullani integral [8, p.472]

(1.13) $$\log u = \int_0^1 \left[\frac{y^{u-1}}{\log y} - \frac{1}{\log y}\right]dy$$

We then have

$$\log u - \psi(u) = \int_0^1 \left[\frac{1}{1-y} + \frac{1}{\log y}\right]y^{u-1}\,dy$$

and, comparing this with (1.8), results in (1.9).

Making the substitution $y = e^{-x}$ in (1.10) we see that

$$I = \int_0^1 \left[\frac{1}{1-y} + \frac{1}{\log y}\right]y^{u-1}\,dy = \int_0^\infty \left[\frac{1}{1-e^{-x}} - \frac{1}{x}\right]e^{-ux}\,dx$$

and noting that

$$\frac{1}{1-e^{-x}} = \frac{e^x}{e^x - 1} = 1 + \frac{1}{e^x - 1}$$

we obtain

$$I = \frac{1}{u} + \int_0^\infty \left[\frac{1}{e^x - 1} - \frac{1}{x}\right]e^{-ux}\,dx$$

Therefore we have from (1.10) [26, p.15]

(1.14) $$\int_0^\infty \left[\frac{1}{e^x - 1} - \frac{1}{x}\right] e^{-ux} dx = -\frac{1}{u} + \log u - \psi(u) = \log u - \psi(1+u)$$

Integrating (1.14) with respect to $u$ over the interval $[0,t]$ gives us a Binet-type integral representation for $\log \Gamma(u)$

(1.15) $$\int_0^\infty \left[\frac{1}{e^x - 1} - \frac{1}{x}\right] \frac{1 - e^{-tx}}{x} dx = t \log t - t - \log \Gamma(1+t)$$

Alternatively, integrating (1.14) with respect to $u$ over the interval $[1,t]$ gives us

(1.16) $$\int_0^\infty \left[\frac{1}{e^x - 1} - \frac{1}{x}\right] \frac{e^{-x} - e^{-tx}}{x} dx = t \log t - t + 1 - \log \Gamma(1+t)$$

which recently appeared in [15].

Further integrations would in turn develop integrals for the Barnes multiple gamma functions defined in [26, p.24]. This is considered further in [12].

Further consideration of integrals such as (1.10) is contained in Appendix B to this paper.

**2. Some integrals involving the Stieltjes constants**

We now boldly assume that (1.6) is also valid in the case where $r$ is a continuous variable (perhaps Carlson's Theorem [27, p.185] could be utilised to provide a rigorous proof: further advice on this aspect would be appreciated). In this regard, it may be noted that Adamchik [1] made a similar assumption in making the transition from equation (23) to equation (24) in his 1997 paper, "A class of logarithmic integrals".

Differentiation of (1.6) with respect to $r$ results in

(2.1) $$\int_0^1 \left[\frac{1}{1-y} + \frac{1}{\log y}\right] \frac{y^{u-1} \log|\log y|}{\log^{r-1} y} dy$$

$$= \frac{1}{\Gamma(r)} \sum_{n=1}^\infty \frac{1}{n+1} \sum_{j=0}^n \binom{n}{j} (-1)^j (u+j)^{r-1} \log^2(u+j) - \frac{\psi(r)}{\Gamma(r)} \sum_{n=1}^\infty \frac{1}{n+1} \sum_{j=0}^n \binom{n}{j} (-1)^j (u+j)^{r-1} \log(u+j)$$

Letting $r = 1$ results in

$$\int_0^1 \left[\frac{1}{1-y} + \frac{1}{\log y}\right] y^{u-1} \log|\log y| dy$$

$$= \sum_{n=1}^{\infty} \frac{1}{n+1} \sum_{j=0}^{n} \binom{n}{j} (-1)^j \log^2(u+j) + \gamma \sum_{n=1}^{\infty} \frac{1}{n+1} \sum_{j=0}^{n} \binom{n}{j} (-1)^j \log(u+j)$$

and starting the summation at $n = 0$ gives us

$$\int_0^1 \left[ \frac{1}{1-y} + \frac{1}{\log y} \right] y^{u-1} \log|\log y| \, dy = \sum_{n=0}^{\infty} \frac{1}{n+1} \sum_{j=0}^{n} \binom{n}{j} (-1)^j \log^2(u+j)$$

$$+ \gamma \sum_{n=0}^{\infty} \frac{1}{n+1} \sum_{j=0}^{n} \binom{n}{j} (-1)^j \log(u+j)) - \log^2 u - \gamma \log u$$

We have previously shown in [13] that for integer $p \geq 0$

(2.2) $$\gamma_p(u) = -\frac{1}{p+1} \sum_{n=0}^{\infty} \frac{1}{n+1} \sum_{j=0}^{n} \binom{n}{j} (-1)^j \log^{p+1}(u+j)$$

where the Stieltjes constants $\gamma_n(u)$ are the coefficients in the Laurent expansion of the Hurwitz zeta function $\varsigma(s,u)$ about $s = 1$

(2.3) $$\varsigma(s,u) = \sum_{n=0}^{\infty} \frac{1}{(n+u)^s} = \frac{1}{s-1} + \sum_{n=0}^{\infty} \frac{(-1)^n}{n!} \gamma_n(u)(s-1)^n$$

and $\gamma_0(u) = -\psi(u)$, where $\psi(u)$ is the digamma function which is the logarithmic derivative of the gamma function $\psi(u) = \frac{d}{du} \log \Gamma(u)$. It is easily seen from the definition of the Hurwitz zeta function that $\varsigma(s,1) = \varsigma(s)$ and accordingly that $\gamma_n(1) = \gamma_n$.

Hence we obtain

(2.4) $$\int_0^1 \left[ \frac{1}{1-y} + \frac{1}{\log y} \right] y^{u-1} \log|\log y| \, dy = -2\gamma_1(u) - \gamma\gamma_0(u) - \log^2 u - \gamma \log u$$

Letting $u = 1$ results in (however see (B.30) in Appendix B)

(2.5) $$\int_0^1 \left[ \frac{1}{1-y} + \frac{1}{\log y} \right] \log|\log y| \, dy = -(2\gamma_1 + \gamma^2)$$

We define $f(y)$ as

$$f(y) = \frac{1}{1-y} + \frac{1}{\log y} = \frac{\log y + 1 - y}{(1-y)\log y}$$

Using L'Hôpital's rule we can easily deduce that $\lim_{y \to 0} f(y) = 1$ and $\lim_{y \to 1} f(y) = -1/2$ and, since $f'(y) < 0$, we note that $f(y)$ is monotonically decreasing on the interval $[0,1]$. However, the range of $\log|\log y|$ is $(-\infty, \infty)$ and hence it is not straightforward to determine the sign of the integral in (2.5). Another approach is therefore necessary.

With the substitution $y = 1/t$ in (2.5) we get

$$\int_0^1 \left[ \frac{1}{1-y} + \frac{1}{\log y} \right] \log|\log y| \, dy = \int_1^\infty \left[ \frac{1}{t-1} - \frac{1}{t \log t} \right] \frac{\log|\log(1/t)|}{t} \, dt$$

We denote $g(t)$ by

$$g(t) = \frac{1}{t-1} - \frac{1}{t \log t} = \frac{t \log t - (t-1)}{t(t-1)\log t} = \frac{\log t - (1 - 1/t)}{(t-1)\log t}$$

and using L'Hôpital's rule we obtain

$$\lim_{t \to 1} g(t) = 1/2 \qquad \lim_{t \to \infty} g(t) = 0$$

Since $\log t > 1 - 1/t$ we see that

$$g(t) > 0 \text{ for } t \in [1, \infty)$$

We see that $h(t) = \dfrac{\log|\log(1/t)|}{t} \to -\infty$ as $t \to 1$ and using L'Hôpital's rule we have

$$\frac{\log|\log(1/t)|}{t} \to 0 \text{ as } t \to \infty$$

We have the derivative

$$h'(t) = \frac{-\dfrac{1}{|\log(1/t)|} - \log|\log(1/t)|}{t^2}$$

Since $\log t > 1 - 1/t$ for $t > 0$ we also have

$$\log|\log(1/t)| > 1 - \frac{1}{|\log(1/t)|} \quad \text{for } |\log(1/t)| > 0$$

and therefore $h'(t)$ is negative with the result that $h(t)$ is monotonic decreasing and negative for $t \in [1,\infty)$. We accordingly deduce that the integrand in (2.5) is negative and that the integral is also negative. This then enables us to conclude that

$$2\gamma_1 + \gamma^2 > 0$$

The constants $\eta_k$ are defined by the following Maclaurin expansion

$$\log[(s-1)\varsigma(s)] = -\sum_{k=0}^{\infty} \frac{\eta_k}{k+1}(s-1)^{k+1} = -\sum_{k=1}^{\infty} \frac{\eta_{k-1}}{k}(s-1)^k$$

where

$$\eta_0 = -\gamma$$

$$\eta_1 = 2\gamma_1 + \gamma^2$$

We note the approximate numerical values of the first three Stieltjes constants from [2]

$$\gamma = 0.5772\cdots \quad \gamma_1 = -0.0728\cdots \quad \gamma_1 = -0.0096\cdots$$

Coffey [9a] has shown that the sequence $(\eta_n)$ has strict sign alteration, i.e.

$$\eta_n = (-1)^{n+1}\varepsilon_n$$

where $\varepsilon_n$ are positive constants. It may be noted that the above analysis has confirmed, without performing any numerical calculations, that $\eta_1$ is indeed positive.

We see that $\dfrac{\log^{2n}|\log(1/t)|}{t} > 0$ and with $j(t) = \dfrac{\log^{2n+1}|\log(1/t)|}{t}$ we see that

$$j'(t) = -\frac{(2n+1)\log^{2n}|\log(1/t)|\dfrac{1}{|\log(1/t)|} + \log^{2n+1}|\log(1/t)|}{t^2}$$

$$= -\frac{\log^{2n}|\log(1/t)|}{t^2}\left(\frac{(2n+1)}{|\log(1/t)|} + \log|\log(1/t)|\right)$$

Since $\log t > 1 - 1/t$ for $t > 0$ we also have

$$\log^{2n+1} t > (1-1/t)^{2n+1} > 1 - \frac{2n+1}{t}$$

where, in the last step, we have employed Bernoulli's inequality. We then have

$$\log^{2n+1} |\log(1/t)| > 1 - \frac{2n+1}{|\log(1/t)|} \quad \text{for } |\log(1/t)| > 0$$

and conclude that $j'(t)$ is negative.

We have therefore determined that

$$(2.6) \quad \int_0^1 \left[\frac{1}{1-y} + \frac{1}{\log y}\right] \log^n |\log y| \, dy = \int_1^\infty \left[\frac{1}{t-1} - \frac{1}{t \log t}\right] \frac{\log^n |\log(1/t)|}{t} \, dt = (-1)^n d_n$$

where $d_n$ are positive constants. This relationship is used in (3.15) below.

□

For convenience, we now write (2.1) as

$$(2.7) \quad \int_0^1 \left[\frac{1}{1-y} + \frac{1}{\log y}\right] \frac{y^{u-1} \log |\log y|}{\log^{r-1} y} \, dy = \frac{1}{\Gamma(r)} S(r,u,2) - \frac{\psi(r)}{\Gamma(r)} S(r,u,1)$$

where $S(r,u,p)$ is defined as

$$S(r,u,p) = \sum_{n=1}^\infty \frac{1}{n+1} \sum_{j=0}^n \binom{n}{j} (-1)^j (u+j)^{r-1} \log^p(u+j)$$

and note that $\frac{d}{dr} S(r,u,p) = S(r,u,p+1)$.

Differentiation of (2.7) with respect to $r$ gives us

$$(2.8) \quad -\int_0^1 \left[\frac{1}{1-y} + \frac{1}{\log y}\right] \frac{y^{u-1} \log^2 |\log y|}{\log^{r-1} y} \, dy = \frac{1}{\Gamma(r)} S(r,u,3) - 2\frac{\psi(r)}{\Gamma(r)} S(r,u,2)$$

$$-\frac{\psi'(r) - [\psi(r)]^2}{\Gamma(r)} S(r,u,1)$$

We additionally define $S_0(r,u,p)$ by

$$S_0(r,u,p) = \sum_{n=0}^{\infty} \frac{1}{n+1} \sum_{j=0}^{n} \binom{n}{j} (-1)^j (u+j)^{r-1} \log^p(u+j)$$

where the summation starts at $n = 0$ and we note that

$$S_0(r,u,p) = S(r,u,p) + u^{r-1} \log^p u$$

and

$$S_0(r,1,p) = S(r,1,p)$$

We also see from (2.2) with $r = 1$ that

$$S(1,u,p+1) = -(p+1)\gamma_p(u) - \log^{p+1} u$$

and

$$S(1,1,p+1) = -(p+1)\gamma_p$$

We then have with $r = 1$ in (2.8)

$$-\int_0^1 \left[ \frac{1}{1-y} + \frac{1}{\log y} \right] y^{u-1} \log^2 |\log y| \, dy = S(1,u,3) - 2\psi(1)S(1,u,2) - \left(\psi'(1) - [\psi(1)]^2\right) S(1,u,1)$$

which is equivalent to

(2.9) $\int_0^1 \left[ \frac{1}{1-y} + \frac{1}{\log y} \right] y^{u-1} \log^2 |\log y| \, dy$

$$= 3\gamma_2(u) + \log^3 u + 2\gamma[2\gamma_1(u) + \log^2 u] - [\varsigma(2) - \gamma^2][\gamma(u) + \log u]$$

Therefore we get with $u = 1$

(2.10) $\int_0^1 \left[ \frac{1}{1-y} + \frac{1}{\log y} \right] \log^2 |\log y| \, dy = 3\gamma_2 + 4\gamma\gamma_1 - [\varsigma(2) - \gamma^2]\gamma$

Therefore we see that

(2.11) $3\gamma_2 + 4\gamma\gamma_1 - [\varsigma(2) - \gamma^2]\gamma > 0$

Coffey [9b] has given the following recurrence relation for the $\eta_n$ constants

$$\eta_n = (-1)^{n+1}\left[\frac{n+1}{n!}\gamma_n + \sum_{k=0}^{n-1}\frac{(-1)^{k+1}}{(n-k-1)!}\gamma_{n-k-1}\eta_k\right]$$

from which we obtain

(2.12) $$\eta_2 = -\left[\frac{3}{2}\gamma_2 - \eta_0\gamma_1 + \eta_1\gamma\right] = -\left[\frac{3}{2}\gamma_2 + 3\gamma\gamma_1 + \gamma^3\right]$$

We then see that

$$\left[\frac{3}{2}\gamma_2 + 3\gamma\gamma_1 + \gamma^3\right] - \gamma\left[\gamma_1 + \frac{1}{2}\gamma + \frac{1}{2}\varsigma(2)\right] > 0$$

but it is not clear if such inequalities are germane in connection with the Riemann Hypothesis via the Li/Keiper criterion.

Differentiation of (2.8) with respect to $r$ gives us

(2.13)
$$\int_0^1\left[\frac{1}{1-y}+\frac{1}{\log y}\right]\frac{y^{u-1}\log^3|\log y|}{\log^{r-1} y}dy = \frac{1}{\Gamma(r)}S(r,u,4) - 3\frac{\psi(r)}{\Gamma(r)}S(r,u,3)$$
$$-3\frac{\psi'(r)-[\psi(r)]^2}{\Gamma(r)}S(r,u,2) - \frac{\psi''(r)-2\psi(r)\psi'(r)+[\psi(r)]^3}{\Gamma(r)}S(r,u,1)$$

Further differentiations will result in integrals $\int_0^1\left[\frac{1}{1-y}+\frac{1}{\log y}\right]y^{u-1}\log^n|\log y|dy$ involving the Stieltjes constants $\gamma_n(u)$ (together with the resulting inequalities).

□

Differentiation of (2.4) with respect to $u$ results in

$$\int_0^1\left[\frac{1}{1-y}+\frac{1}{\log y}\right]y^{u-1}\log y\cdot\log|\log y|dy = -2\gamma_1'(u) - \gamma\gamma_0'(u) - \frac{\gamma+2\log u}{u}$$

and with $u=1$ we get using $\gamma_0'(u)=\psi'(u)$

$$-\int_0^1\left[\frac{1}{1-y}+\frac{1}{\log y}\right]\log y\cdot\log|\log y|dy = 2\gamma_1'(1) - \gamma\varsigma(2) + \gamma$$

Using equation (4.3.244) in [13]

$$\gamma_1'(1) = 2\pi^2 \varsigma'(-1) + \varsigma(2)(\gamma + \log 2\pi)$$

we may write this as

(2.14) $\quad -\int_0^1 \left[\frac{1}{1-y} + \frac{1}{\log y}\right] \log y \cdot \log|\log y| \, dy = 4\pi^2 \varsigma'(-1) + \varsigma(2)(\gamma + 2\log 2\pi) + \gamma$

Using (C.1) this may be written as

(2.15) $\quad -\int_0^1 \frac{\log y \cdot \log|\log y|}{1-y} \, dy = 4\pi^2 \varsigma'(-1) + \varsigma(2)(\gamma + 2\log 2\pi)$

## 3. An application of the (exponential) complete Bell polynomials

Alternatively, we may write (1.6) as

(3.1) $\quad \Gamma(r) \int_0^1 \left[\frac{1}{1-y} + \frac{1}{\log y}\right] \frac{y^{u-1}}{\log^{r-1} y} \, dy = -\sum_{n=1}^{\infty} \frac{1}{n+1} \sum_{j=0}^{n} \binom{n}{j} (-1)^j (u+j)^{r-1} \log(u+j)$

Letting $h(r)$ be defined as

$$h(r) = \frac{\Gamma(r)}{\log^{r-1} y}$$

we have

$$\log h(r) = \log \Gamma(r) - (r-1) \log|\log y|$$

and differentiation results in

(3.2) $\quad h'(r) = h(r)[\psi(r) - \log|\log y|]$

Using (3.2) we now differentiate (3.1) with respect to $r$ to obtain

$$\Gamma(r) \int_0^1 \left[\frac{1}{1-y} + \frac{1}{\log y}\right] \frac{y^{u-1}[\psi(r) - \log|\log y|]}{\log^{r-1} y} \, dy = -\sum_{n=1}^{\infty} \frac{1}{n+1} \sum_{j=0}^{n} \binom{n}{j} (-1)^j (u+j)^{r-1} \log^2(u+j)$$

which may be expressed as

$$\Gamma(r)\int_0^1 \left[\frac{1}{1-y}+\frac{1}{\log y}\right]\frac{y^{u-1}[\psi(r)-\log|\log y|]}{\log^{r-1} y}dy = 2\gamma_1(u)+u^{r-1}\log^2 u$$

Letting $r=1$ this becomes

$$\int_0^1 \left[\frac{1}{1-y}+\frac{1}{\log y}\right]y^{u-1}[\gamma+\log|\log y|]dy = -2\gamma_1(u)-\log^2 u$$

and with $u=1$ we have

$$\int_0^1 \left[\frac{1}{1-y}+\frac{1}{\log y}\right][\gamma+\log|\log y|]dy = -2\gamma_1$$

which of course is consistent with our previous results.

We now wish to obtain a formula for the $m$ th derivative of (3.1).

It is shown in Appendix A that

(3.3) $$\frac{d^m}{dx^m}e^{f(x)} = e^{f(x)}Y_m\left(f^{(1)}(x), f^{(2)}(x),..., f^{(m)}(x)\right)$$

where the (exponential) complete Bell polynomials $Y_n(x_1,...,x_n)$ are defined by $Y_0=1$ and for $n\geq 1$

(3.4) $$Y_n(x_1,...,x_n) = \sum_{\pi(n)} \frac{n!}{k_1!\, k_2!...\, k_n!}\left(\frac{x_1}{1!}\right)^{k_1}\left(\frac{x_2}{2!}\right)^{k_2}...\left(\frac{x_n}{n!}\right)^{k_n}$$

where the sum is taken over all partitions $\pi(n)$ of $n$, i.e. over all sets of integers $k_j$ such that

$$k_1+2k_2+3k_3+...+nk_n = n$$

Suppose that $h'(x)=h(x)g(x)$ and let $f(x)=\log h(x)$. We see that

$$f'(x) = \frac{h'(x)}{h(x)} = g(x)$$

and then using (3.3) above we have

(3.5) $$\frac{d^m}{dx^m}h(x) = \frac{d^m}{dx^m}e^{\log h(x)} = h(x)Y_m\left(g(x), g^{(1)}(x),..., g^{(m-1)}(x)\right)$$

As an example, letting $h(x) = \Gamma(x)$ in (3.5) we obtain

(3.6) $$\frac{d^m}{dx^m}e^{\log \Gamma(x)} = \Gamma^{(m)}(x) = \Gamma(x)Y_m\left(\psi(x), \psi^{(1)}(x),..., \psi^{(m-1)}(x)\right)$$

$$= \int_0^\infty t^{x-1}e^{-t}\log^m t \, dt$$

and since [26, p.22]

(3.7) $$\psi^{(p)}(x) = (-1)^{p+1}p!\varsigma(p+1, x)$$

we may express $\Gamma^{(m)}(x)$ in terms of $\psi(x)$ and the Hurwitz zeta functions. In particular, Kölbig [20] noted that

(3.8) $$\Gamma^{(m)}(1) = Y_m(-\gamma, x_1,..., x_{m-1})$$

where $x_p = (-1)^{p+1}p!\varsigma(p+1)$. Values of $\Gamma^{(m)}(1)$ for $m \leq 10$ are reported in [26, p.265] and the first three are

$$\Gamma^{(1)}(1) = -\gamma$$

$$\Gamma^{(2)}(1) = \varsigma(2) + \gamma^2$$

$$\Gamma^{(3)}(1) = -[2\varsigma(3) + 3\gamma\varsigma(2) + \gamma^3]$$

As shown in Appendix C, we note that $\Gamma^{(n)}(1)$ has the same sign as $(-1)^n$. This was also reported as an exercise in Apostol's book [4, p.303].

As a variation of (3.5) above, suppose that $j'(x) = j(x)[g(x) + \alpha]$ where $\alpha$ is independent of $x$ and let $f(x) = \log j(x)$. We see that

$$f'(x) = \frac{j'(x)}{j(x)} = g(x) + \alpha \text{ and } f^{(k+1)}(x) = g^{(k+1)}(x) \text{ for } k \geq 1$$

and therefore we obtain

(3.9) $$\frac{d^m}{dx^m} j(x) = \frac{d^m}{dx^m} e^{\log j(x)} = j(x) Y_m\left(g(x)+\alpha, g^{(1)}(x),...,g^{(m-1)}(x)\right)$$

It is also shown in Appendix A that

$$Y_n(x_1+\alpha,...,x_n) = \sum_{k=0}^{n}\binom{n}{k}\alpha^{n-k} Y_k(x_1,...,x_k)$$

and we then determine that

(3.10) $$\frac{d^m}{dx^m} j(x) = j(x)\sum_{k=0}^{m}\binom{m}{k}\alpha^{m-k} Y_k\left(g(x), g^{(1)}(x),...,g^{(k-1)}(x)\right)$$

Now, referring back to (3.2), we see that with $g(r) = \psi(r) - \log|\log y|$

$$\frac{d^m}{dr^m}\frac{\Gamma(r)}{\log^{r-1} y} = \frac{\Gamma(r)}{\log^{r-1} y}\sum_{k=0}^{m}\binom{m}{k}(-1)^{m-k}\log^{m-k}|\log y|\cdot Y_k\left(\psi(r),\psi^{(1)}(r),...,\psi^{(k-1)}(r)\right)$$

and since we have from (3.5)

$$\Gamma^{(k)}(r) = \Gamma(r)Y_k\left(\psi(r),\psi^{(1)}(r),...,\psi^{(k-1)}(r)\right)$$

we obtain

$$\frac{d^m}{dr^m}\frac{\Gamma(r)}{\log^{r-1} y} = \frac{\Gamma(r)}{\log^{r-1} y}\sum_{k=0}^{m}\binom{m}{k}(-1)^{m-k}\Gamma^{(k)}(r)\log^{m-k}|\log y|$$

Hence, differentiating (3.1) $m$ times with respect to $r$ results in

(3.11) $$\Gamma(r)\sum_{k=0}^{m}\binom{m}{k}(-1)^{m-k}\Gamma^{(k)}(r)\int_0^1\left[\frac{1}{1-y}+\frac{1}{\log y}\right]\frac{y^{u-1}\log^{m-k}|\log y|}{\log^{r-1} y}dy$$

$$= -\sum_{n=1}^{\infty}\frac{1}{n+1}\sum_{j=0}^{n}\binom{n}{j}(-1)^j(u+j)^{r-1}\log^{m+1}(u+j)$$

Letting $r=1$ this becomes

$$\sum_{k=0}^{m}\binom{m}{k}(-1)^{m-k}\Gamma^{(k)}(1)\int_0^1\left[\frac{1}{1-y}+\frac{1}{\log y}\right]y^{u-1}\log^{m-k}|\log y|dy$$

$$= -\sum_{n=1}^{\infty} \frac{1}{n+1} \sum_{j=0}^{n} \binom{n}{j} (-1)^j \log^{m+1}(u+j)$$

$$= (m+1)\gamma_m(u) + \log^{m+1} u$$

With $u = 1$ we obtain

(3.12) $\quad (m+1)\gamma_m = \sum_{k=0}^{m} \binom{m}{k} (-1)^{m-k} \Gamma^{(k)}(1) \int_0^1 \left[ \frac{1}{1-y} + \frac{1}{\log y} \right] \log^{m-k} |\log y| \, dy$

For example, with $m = 1$ we get

$$2\gamma_1 = -\int_0^1 \left[ \frac{1}{1-y} + \frac{1}{\log y} \right] \log |\log y| \, dy - \gamma \int_0^1 \left[ \frac{1}{1-y} + \frac{1}{\log y} \right] dy$$

which is in agreement with our previous results.

Equation (3.12) sheds a little light upon the complexity involved in the alteration in the signs of $\gamma_m$ because

(3.13) $\quad \Gamma^{(k)}(1) = (-1)^k c_k$

(3.14) $\quad \int_0^1 \left[ \frac{1}{1-y} + \frac{1}{\log y} \right] \log^{m-k} |\log y| \, dy = (-1)^{m-k} d_{m-k}$

where $c_k$ and $d_k$ are positive constants. The relationship (3.13) is derived in Appendix C and (3.14) was shown in (2.6) above. We then have

(3.15) $\quad (m+1)\gamma_m = \sum_{k=0}^{m} \binom{m}{k} (-1)^k c_k d_{m-k}$

$\square$

Referring back to (1.1) we make the summation

$$\sum_{n=0}^{\infty} \frac{1}{n+1} \int_0^{\infty} e^{-ax} (1-e^{-x})^n x^{p-1} dx = \Gamma(p) \sum_{n=0}^{\infty} \frac{1}{n+1} \sum_{k=0}^{n} \binom{n}{k} \frac{(-1)^k}{(a+k)^p}$$

and we have

$$\sum_{n=0}^{\infty}\frac{1}{n+1}\int_0^{\infty} e^{-ax}(1-e^{-x})^n x^{p-1}dx = \int_0^{\infty} e^{-ax} x^{p-1} \sum_{n=0}^{\infty}\frac{(1-e^{-x})^n}{n+1} dx$$

$$= \int_0^{\infty}\frac{e^{-ax} x^{p-1}}{1-e^{-x}} \sum_{n=0}^{\infty}\frac{(1-e^{-x})^{n+1}}{n+1} dx$$

$$= \int_0^{\infty}\frac{e^{-ax} x^{p-1} \log(1-e^{-x})}{1-e^{-x}} dx$$

Therefore we obtain

$$\int_0^{\infty}\frac{e^{-ax} x^{p-1} \log(1-e^{-x})}{1-e^{-x}} dx = \Gamma(p) \sum_{n=0}^{\infty}\frac{1}{n+1}\sum_{k=0}^{n}\binom{n}{k}\frac{(-1)^k}{(a+k)^p}$$

and with $a=1$ we have

$$\int_0^{\infty}\frac{e^{-x} x^{p-1} \log(1-e^{-x})}{1-e^{-x}} dx = \Gamma(p) \sum_{n=0}^{\infty}\frac{1}{n+1}\sum_{k=0}^{n}\binom{n}{k}\frac{(-1)^k}{(1+k)^p}$$

We note that the Bernoulli polynomials may be represented by [17]

$$B_k(a) = \sum_{n=0}^{\infty}\frac{1}{n+1}\sum_{k=0}^{n}\binom{n}{k}(-1)^k (a+k)^k$$

and hence for $m \geq 1$

$$B_{2m+1}(1) = \sum_{n=0}^{\infty}\frac{1}{n+1}\sum_{k=0}^{n}\binom{n}{k}(-1)^k (1+k)^{2m+1} = 0$$

We now consider the limit as $p \to -(2m+1)$ where $m$ is a positive integer and $m \geq 1$; the right-hand side of () is indeterminate but we may use L'Hôpital's rule when considering the ratio

$$\frac{\sum_{n=0}^{\infty}\frac{1}{n+1}\sum_{k=0}^{n}\binom{n}{k}\frac{(-1)^k}{(a+k)^p}}{1/\Gamma(p)}$$

As shown previously, using Euler's reflection formula we have

$$\frac{d}{dp}[1/\Gamma(p)] = \Gamma(1-p)\cos \pi p - \frac{1}{\pi}\Gamma'(1-p)\sin \pi p$$

and therefore we see that where $p \to -(2m+1)$

$$\lim_{p \to -(2m+1)} \frac{d}{dp}[1/\Gamma(p)] = (-1)^{2m+1}\Gamma(2m)$$

We thereby obtain

$$\int_0^\infty \frac{e^{-x}\log(1-e^{-x})}{(1-e^{-x})x^{2m}}dx = (-1)^{2m+1}\Gamma(2m)\sum_{n=0}^\infty \frac{1}{n+1}\sum_{k=0}^n \binom{n}{k}(-1)^k(1+k)^{2m+1}\log(1+k)$$

# Appendix A

## A brief survey of the (exponential) complete Bell polynomials

The (exponential) complete Bell polynomials are defined by $Y_0 = 1$ and for $n \geq 1$

(A1) $\qquad Y_n(x_1,...,x_n) = \sum_{\pi(n)} \frac{n!}{k_1! k_2!... k_n!}\left(\frac{x_1}{1!}\right)^{k_1}\left(\frac{x_2}{2!}\right)^{k_2}...\left(\frac{x_n}{n!}\right)^{k_n}$

where the sum is taken over all partitions $\pi(n)$ of $n$, i.e. over all sets of integers $k_j$ such that

$$k_1 + 2k_2 + 3k_3 + ... + nk_n = n$$

The complete Bell polynomials have integer coefficients and the first six are set out below [11, p.307]

(A.2) $\qquad Y_1(x_1) = x_1$

$$Y_2(x_1, x_2) = x_1^2 + x_2$$

$$Y_3(x_1, x_2, x_3) = x_1^3 + 3x_1 x_2 + x_3$$

$$Y_4(x_1, x_2, x_3, x_4) = x_1^4 + 6x_1^2 x_2 + 4x_1 x_3 + 3x_2^2 + x_4$$

$$Y_5(x_1, x_2, x_3, x_4, x_5) = x_1^5 + 10x_1^3 x_2 + 10x_1^2 x_3 + 15x_1 x_2^2 + 5x_1 x_4 + 10x_2 x_3 + x_5$$

$$Y_6(x_1, x_2, x_3, x_4, x_5, x_6) = x_1^6 + 6x_1 x_5 + 15x_2 x_4 + 10x_3^2 + 15x_1^2 x_4 + 15x_2^3 + 60x_1 x_2 x_3$$

$$+20x_1^3 x_3 + 45 x_1^2 x_2^2 + 15 x_1^4 x_1 + x_6$$

The total number of terms $\pi(n)$ increases rapidly; for example, as reported by Bell [5] in 1934, we have $\pi(22) = 1002$ terms.

The modus operandi of the summation in (A.1) is easily illustrated by the following example: if $n = 4$, then $k_1 + 2k_2 + 3k_3 + 4k_4 = 4$ is satisfied by the integers in the following array

$$\begin{bmatrix} k_1 & k_2 & k_3 & k_4 \\ 4 & 0 & 0 & 0 \\ 2 & 1 & 0 & 0 \\ 1 & 0 & 1 & 0 \\ 0 & 2 & 0 & 0 \\ 0 & 0 & 0 & 1 \end{bmatrix}$$

The complete Bell polynomials are also given by the exponential generating function (Comtet [11, p.134])

(A.3) $$\exp\left(\sum_{j=1}^{\infty} x_j \frac{t^j}{j!}\right) = 1 + \sum_{n=1}^{\infty} Y_n(x_1, ..., x_n) \frac{t^n}{n!} = \sum_{n=0}^{\infty} Y_n(x_1, ..., x_n) \frac{t^n}{n!}$$

Let us now consider a function $f(x)$ which has a Taylor series expansion around $x$: we have

$$e^{f(x+t)} = \exp\left(\sum_{j=0}^{\infty} f^{(j)}(x) \frac{t^j}{j!}\right) = e^{f(x)} \exp\left(\sum_{j=1}^{\infty} f^{(j)}(x) \frac{t^j}{j!}\right)$$

$$= e^{f(x)} \left\{ 1 + \sum_{n=1}^{\infty} Y_n\left(f^{(1)}(x), f^{(2)}(x), ..., f^{(n)}(x)\right) \frac{t^n}{n!} \right\}$$

We see that

$$\frac{d^m}{dx^m} e^{f(x)} = \frac{\partial^m}{\partial x^m} e^{f(x+t)} \bigg|_{t=0} = \frac{\partial^m}{\partial t^m} e^{f(x+t)} \bigg|_{t=0}$$

and we therefore obtain (as noted by Kölbig [20])

(A.4) $$\frac{d^m}{dx^m} e^{f(x)} = e^{f(x)} Y_m\left(f^{(1)}(x), f^{(2)}(x), \ldots, f^{(m)}(x)\right)$$

Differentiating (A.4) we see that

$$Y_{m+1}\left(f^{(1)}(x), f^{(2)}(x), \ldots, f^{(m+1)}(x)\right) = \left(f^{(1)}(x) + \frac{d}{dx}\right) Y_m\left(f^{(1)}(x), f^{(2)}(x), \ldots, f^{(m)}(x)\right)$$

Suppose that $h'(x) = h(x)g(x)$ and let $f(x) = \log h(x)$. We see that

$$f'(x) = \frac{h'(x)}{h(x)} = g(x)$$

and then using (A.4) above we have

(A.5) $$\frac{d^m}{dx^m} h(x) = \frac{d^m}{dx^m} e^{\log h(x)} = h(x) Y_m\left(g(x), g^{(1)}(x), \ldots, g^{(m-1)}(x)\right)$$

As an example, letting $f(x) = \log \Gamma(x)$ in (A.4) we obtain

(A.6) $$\frac{d^m}{dx^m} e^{\log \Gamma(x)} = \Gamma^{(m)}(x) = \Gamma(x) Y_m\left(\psi(x), \psi^{(1)}(x), \ldots, \psi^{(m-1)}(x)\right)$$

$$= \int_0^\infty t^{x-1} e^{-t} \log^m t \, dt$$

and since [26, p.22]

(A.7) $$\psi^{(p)}(x) = (-1)^{p+1} p! \varsigma(p+1, x)$$

we may express $\Gamma^{(m)}(x)$ in terms of $\psi(x)$ and the Hurwitz zeta functions. In particular, Kölbig [20] notes that

(A.8) $$\Gamma^{(m)}(1) = Y_m(-\gamma, x_1, \ldots, x_{m-1})$$

where $x_p = (-1)^{p+1} p! \varsigma(p+1)$. Values of $\Gamma^{(m)}(1)$ are reported in [26, p.265] for $m \leq 10$ and the first three are

$$\Gamma^{(1)}(1) = -\gamma$$

$$\Gamma^{(2)}(1) = \varsigma(2) + \gamma^2$$

$$\Gamma^{(3)}(1) = -[2\varsigma(3) + 3\gamma\varsigma(2) + \gamma^3]$$

The general form is

$$\Gamma^{(m)}(1) = (-1)^m \sum_{j=1}^{m} \varepsilon_{mj}$$

where $\varepsilon_{mj}$ are positive constants.

We could also, for example, let $f(x) = \log \sin(\pi x)$ in (A.4) to obtain complete Bell polynomials involving the derivatives of $\cot(\pi x)$.

We now consider a minor modification of the arguments of the exponential generating function

$$\exp\left(\sum_{j=1}^{\infty} x_j \frac{t^j}{j!}\right) = \sum_{n=0}^{\infty} Y_n(x_1,...,x_n) \frac{t^n}{n!}$$

where we let $x_1 \to = x_1 + \alpha$ and the other arguments remain unchanged. We then see that

$$\exp\left(\alpha t + \sum_{j=1}^{\infty} x_j \frac{t^j}{j!}\right) = \sum_{n=0}^{\infty} Y_n(x_1 + \alpha,...,x_n) \frac{t^n}{n!}$$

Since

$$\exp\left(\alpha t + \sum_{j=1}^{\infty} x_j \frac{t^j}{j!}\right) = \exp(\alpha t) \exp\left(\sum_{j=1}^{\infty} x_j \frac{t^j}{j!}\right)$$

we easily determine that

(A.9) $$\exp(\alpha t) \sum_{n=0}^{\infty} Y_n(x_1,...,x_n) \frac{t^n}{n!} = \sum_{n=0}^{\infty} Y_n(x_1 + \alpha,...,x_n) \frac{t^n}{n!}$$

Using the Cauchy series product formula we find that

$$\sum_{n=0}^{\infty} \alpha^n \frac{t^n}{n!} \sum_{n=0}^{\infty} Y_n(x_1,...,x_n) \frac{t^n}{n!} = \sum_{n=0}^{\infty} \sum_{k=0}^{n} \frac{\alpha^{n-k}}{(n-k)!} \frac{Y_k(x_1,...,x_k)}{k!} t^n$$

and equating coefficients we see that

(A.10) $$Y_n(x_1 + \alpha,...,x_n) = \sum_{k=0}^{n} \binom{n}{k} \alpha^{n-k} Y_k(x_1,...,x_k)$$

□

As a variation of (A.5) above, suppose that $j'(x) = j(x)[g(x)+\alpha]$ where $\alpha$ is independent of $x$ and let $f(x) = \log j(x)$. We see that

$$f'(x) = \frac{j'(x)}{j(x)} = g(x) + \alpha \text{ and } f^{(k+1)}(x) = g^{(k)}(x) \text{ for } k \geq 1$$

and therefore we obtain

(A.11) $$\frac{d^m}{dx^m} j(x) = \frac{d^m}{dx^m} e^{\log j(x)} = j(x) Y_m\left(g(x)+\alpha, g^{(1)}(x),...,g^{(m-1)}(x)\right)$$

Using (A.10) we then determine that

(A.12) $$\frac{d^m}{dx^m} j(x) = j(x) \sum_{k=0}^{m} \binom{m}{k} \alpha^{m-k} Y_k\left(g(x), g^{(1)}(x),...,g^{(k-1)}(x)\right)$$

The relation (A.10) may be generalised as follows to

(A.13) $$Y_n(x_1+y_1,...,x_n+y_n) = \sum_{k=0}^{n} \binom{n}{k} Y_{n-k}(x_1,...,x_{n-k}) Y_k(y_1,...,y_k)$$

and we note that

$$Y_k(\alpha,0,...,0) = \alpha^k$$

Equation (A.13) follows by noting that

$$Y_n(x_1+y_1,...,x_n+y_n) = \exp\left(\sum_{j=1}^{\infty}(x_j+y_j)\frac{t^j}{j!}\right) = \exp\left(\sum_{j=1}^{\infty} x_j \frac{t^j}{j!}\right) \exp\left(\sum_{j=1}^{\infty} y_j \frac{t^j}{j!}\right)$$

$$= \sum_{n=0}^{\infty} Y_n(x_1,...,x_n)\frac{t^n}{n!} \cdot \sum_{n=0}^{\infty} Y_n(y_1,...,y_n)\frac{t^n}{n!}$$

and, as before, we apply the Cauchy series product formula.

# Appendix B

# Miscellaneous integral identities

In [24] Rivoal recently gave an elementary proof of the following lemma:

(B.1) $$\int_0^1 x^n \Omega(x)\,dx = \gamma - \left[H_n - \log(n+1)\right]$$

where, for convenience $\Omega(x)$, is defined as $\Omega(x) = \dfrac{1}{1-x} + \dfrac{1}{\log x}$ and $H_n$ is the harmonic number defined by

$$H_n = \sum_{k=1}^{n} \frac{1}{k}$$

Using L'Hôpital's rule it is easily seen that $\Omega(x)$ is continuous on $[0,1]$. For $x \in (0,1)$ we also have by a straightforward integration

$$\Omega(x) = \int_0^1 \frac{1-x^t}{1-x}\,dt$$

We see that

$$\int_0^1 \Omega(x)\,dx = \int_0^1 \left(\frac{1}{1-x} + \frac{1}{\log x}\right)dx = \int_0^\infty \left[\frac{1}{1-e^{-u}} - \frac{1}{u}\right] e^{-u}\,du = \gamma$$

Using the well-known integral representation of the Riemann zeta function

$$\varsigma(s)\Gamma(s) = \int_0^\infty \frac{u^{s-1}}{e^u - 1}\,du$$

we also see that $\Gamma(s) = (s-1)\Gamma(s-1) = (s-1)\int_0^\infty e^{-u} u^{s-2}\,du$. Therefore we get

$$\left[\varsigma(s) - \frac{1}{s-1}\right]\Gamma(s) = \int_0^\infty u^{s-1}\left[\frac{1}{e^u-1} - \frac{1}{ue^u}\right]du$$

Hence, in the limit as $s \to 1$, we obtain the familiar limit

$$\lim_{s \to 1}\left[\varsigma(s) - \frac{1}{s-1}\right]\Gamma(s) = \lim_{s \to 1}\left[\varsigma(s) - \frac{1}{s-1}\right] = \int_0^\infty \left[\frac{1}{1-e^{-u}} - \frac{1}{u}\right]e^{-u}\,du = \gamma$$

We now continue with Rivoal's lemma.

We have by simple algebra

$$\int_0^1 x^n \Omega(x)\,dx = \int_0^1 \Omega(x)\,dx - \int_0^1 \frac{x^n-1}{x-1}\,dx + \int_0^1 \frac{x^n-1}{\log x}\,dx$$

and this gives us

(B.2) $\quad \int_0^1 x^n\left(\frac{1}{1-x}+\frac{1}{\log x}\right)dx = \gamma - H_n + \int_0^1 \frac{x^n-1}{\log x}\,dx = \gamma - [H_n - \log(n+1)]$

It may be noted that this is equivalent to letting $u = n+1$ in (1.10) because

$$\psi(n+1) = -\gamma + H_n$$

Hence by summation of (B.2) we obtain for $\operatorname{Re}(s) > 1$

(B.3) $\quad \gamma\varsigma(s) - \sum_{n=1}^{\infty}\frac{H_n}{n^s} + \sum_{n=1}^{\infty}\frac{\log(n+1)}{n^s} = \int_0^1 Li_s(x)\left(\frac{1}{1-x}+\frac{1}{\log x}\right)dx$

where $Li_s(x)$ is the polylogarithm function defined by

$$Li_s(x) = \sum_{n=1}^{\infty}\frac{x^n}{n^s}$$

In 1997 Candelpergher et al. [9] produced a somewhat similar result in the case where $s = 2$ by reference to "Ramanujan summation" but unfortunately I am not au fait with the underlying analysis contained in that paper.

☐

We have from [23]

(B.4) $\quad H_n - \log n - \gamma = \int_0^{\infty} e^{-nx}\left[\frac{1}{x} - \frac{1}{e^x-1}\right]dx$

and on summation we obtain

$$\sum_{n=1}^{\infty}\frac{H_n}{n^s} + \varsigma'(s) - \gamma\varsigma(s) = \sum_{n=1}^{\infty}\frac{1}{n^s}\int_0^{\infty} e^{-nx}\left[\frac{1}{x} - \frac{1}{e^x-1}\right]dx$$

In particular we have

$$\sum_{n=1}^{\infty}\frac{1}{n^2}\int_0^{\infty} e^{-nx}\left[\frac{1}{x}-\frac{1}{e^x-1}\right]dx = \sum_{n=1}^{\infty}\int_0^{\infty} e^{-nu}du\int_0^{\infty} e^{-nv}dv\int_0^{\infty} e^{-nx}\left[\frac{1}{x}-\frac{1}{e^x-1}\right]dx$$

$$= \sum_{n=1}^{\infty}\int_0^{\infty}du\int_0^{\infty}dv\int_0^{\infty} e^{-n(x+u+v)}\left[\frac{1}{x}-\frac{1}{e^x-1}\right]dx$$

$$= \int_0^{\infty}du\int_0^{\infty}dv\int_0^{\infty}\frac{1}{e^{x+u+v}-1}\left[\frac{1}{x}-\frac{1}{e^x-1}\right]dx$$

We have

$$\int_0^{\infty}\frac{1}{e^{x+u+v}-1}du = \int_0^{\infty}\frac{e^{-(x+u+v)}}{1-e^{-(x+u+v)}}du = \log[1-e^{-(x+u+v)}]\Big|_0^{\infty} = -\log[1-e^{-(x+v)}]$$

and the Wolfram Integrator tells us that

$$\int \log[1-e^{-(x+v)}]dv = \frac{1}{2}v^2 + v\log\frac{[1-e^{-(x+v)}]}{[1-e^{(x+v)}]} - Li_2[e^{(x+v)}]$$

With the substitution $y = e^{-v}$ we can easily find that

$$\int_0^{\infty}\log[1-e^{-(x+v)}]dv = -Li_2[e^{-x}]$$

and hence we obtain

$$\sum_{n=1}^{\infty}\frac{H_n}{n^2} + \varsigma'(2) - \gamma\varsigma(2) = \int_0^{\infty} Li_2[e^{-x}]\left[\frac{1}{x}-\frac{1}{e^x-1}\right]dx$$

or equivalently

(B.5) $$2\varsigma(3) + \varsigma'(2) - \gamma\varsigma(2) = \int_0^{\infty} Li_2[e^{-x}]\left[\frac{1}{x}-\frac{1}{e^x-1}\right]dx$$

where we have employed the well-known Euler sum [26, p.103] $\sum_{n=1}^{\infty}\frac{H_n}{n^2} = 2\varsigma(3)$.

With the substitution $t = e^{-x}$ this may be written as (B.6) $\quad 2\varsigma(3)+\varsigma'(2)-\gamma\varsigma(2) = \int_0^1 Li_2(t)\left[\frac{1}{t-1}-\frac{1}{t\log t}\right]dt$

We see from (1.13) that

$$\sum_{n=1}^\infty \frac{\log n}{n^s} = \sum_{n=1}^\infty \frac{1}{n^s}\int_0^\infty \frac{e^{-x}-e^{-nx}}{x}dx$$

which may be written as

$$-\varsigma'(s) = \int_0^\infty \frac{\varsigma(s)e^{-x}-Li_s[e^{-x}]}{x}dx$$

With $t = e^{-x}$ this becomes

(B.7) $\quad\varsigma'(s) = \int_0^1 \frac{\varsigma(s)t - Li_s(t)}{t\log t}dt$

and with $s = 2$ we get

(B.8) $\quad\varsigma'(2) = \int_0^1 \frac{\varsigma(2)t - Li_2(t)}{t\log t}dt$

We note from (B.6) that

(B.9) $\quad 2\varsigma(3)+\varsigma'(2)-\gamma\varsigma(2) = \int_0^1 Li_2(t)\left[\frac{1}{t-1}-\frac{1}{t\log t}\right]dt$

$$= \int_0^1 \left[\frac{Li_2(t)}{t-1} - \frac{Li_2(t)}{t\log t}\right]dt$$

$$= \int_0^1 \left[\frac{Li_2(t)}{t-1} - \frac{\varsigma(2)}{\log t} + \frac{\varsigma(2)}{\log t} - \frac{Li_2(t)}{t\log t}\right]dt$$

$$= \int_0^1 \left[\frac{Li_2(t)}{t-1} - \frac{\varsigma(2)}{\log t}\right]dt + \int_0^1 \left[\frac{\varsigma(2)}{\log t} - \frac{Li_2(t)}{t\log t}\right]dt$$

and using (B.8) this becomes

$$= \int_0^1 \left[ \frac{Li_2(t)}{t-1} - \frac{\varsigma(2)}{\log t} \right] dt + \varsigma'(2)$$

We therefore see that

(B.10) $$2\varsigma(3) - \gamma\varsigma(2) = \int_0^1 \left[ \frac{Li_2(t)}{t-1} - \frac{\varsigma(2)}{\log t} \right] dt$$

We note from (B.3) with $s = 2$ that

$$2\varsigma(3) - \sum_{n=1}^{\infty} \frac{\log(n+1)}{n^2} - \gamma\varsigma(2) = \int_0^1 Li_2(t) \left( \frac{1}{t-1} - \frac{1}{\log t} \right) dt$$

Subtracting (B.6) results in

(B.11) $$\varsigma'(2) + \sum_{n=1}^{\infty} \frac{\log(n+1)}{n^2} = \int_0^1 \frac{(t-1)Li_2(t)}{t \log t} dt$$

or equivalently

$$\sum_{n=1}^{\infty} \frac{\log(n+1)}{n^2} - \sum_{n=1}^{\infty} \frac{\log n}{n^2} = \int_0^1 \frac{(t-1)Li_2(t)}{t \log t} dt$$

$\square$

We may write (B.10) as the limit

$$2\varsigma(3) - \gamma\varsigma(2) = \lim_{x \to 1} \int_0^x \left[ \frac{Li_2(t)}{t-1} - \frac{\varsigma(2)}{\log t} \right] dt$$

Using integration by parts we can easily determine that

$$\int \frac{Li_2(t)}{t-1} dt = Li_2(t)\log(1-t) + \log t \log^2(1-t) + 2Li_2(1-t)\log(1-t) - 2Li_3(1-t)$$

and we therefore have

$$\int_0^x \frac{Li_2(t)}{t-1} dt = Li_2(x)\log(1-x) + \log x \log^2(1-x) + 2Li_2(1-x)\log(1-x) - 2Li_3(1-x) + 2\varsigma(3)$$

Using L'Hôpital's rule we find that

$$\lim_{x \to 1}\left[(1-x)\log(1-x)\right] = 0$$

and hence we have

$$\lim_{x \to 1}\left[Li_2(1-x)\log(1-x)\right] = 0$$

We therefore have the limit

$$\lim_{x \to 1} \int_0^x \frac{Li_2(t)}{t-1} dt = \lim_{x \to 1}[Li_2(x)\log(1-x) + \log x \log^2(1-x)] + 2\varsigma(3)$$

which then implies that

$$-\gamma\varsigma(2) = \lim_{x \to 1}\left[Li_2(x)\log(1-x) + \log x \log^2(1-x) - \int_0^x \frac{\varsigma(2)}{\log t} dt\right]$$

We therefore have

(B.12) $$\gamma\varsigma(2) = \lim_{x \to 1}\left[-Li_2(x)\log(1-x) - \log x \log^2(1-x) + \varsigma(2)li(x)\right]$$

where $li(x)$ is the logarithmic integral. Nielsen [21, p.3] has shown that for $0 < x < 1$

(B.13) $$li(x) = \gamma + \log(-\log x) + \sum_{n=1}^{\infty} \frac{\log^n x}{n!n}$$

and hence we have

$$\gamma\varsigma(2) = \lim_{x \to 1}\left[-Li_2(x)\log(1-x) - \log x \log^2(1-x) + \varsigma(2)[\gamma + \log(-\log x)]\right]$$

which gives us

$$\lim_{x \to 1}\left[-Li_2(x)\log(1-x) - \log x \log^2(1-x) + \varsigma(2)\log(-\log x)\right] = 0$$

Using Euler's identity

$$\varsigma(2) = \log x \log(1-x) + Li_2(x) + Li_2(1-x)$$

we may write this as

$$\lim_{x \to 1}\left[Li_2(1-x)\log(1-x)+\varsigma(2)\log(1-x)+\varsigma(2)\log(-\log x)\right]=0$$

and hence we have

(B.14) $$\lim_{x \to 1}\left[\log(1-x)+\log(-\log x)\right]=0$$

Since

$$\frac{d}{dt}\log[\log(1/t)]=-\frac{1}{t\log(1/t)}=\frac{1}{t\log t}$$

integration by parts gives us

$$\int_0^x \log[\log(1/t)]dt = t\log[\log(1/t)]\Big|_0^x - \int_0^x \frac{dt}{\log t}$$

Using L'Hôpital's rule we see that

$$\lim_{t \to 0} t\log[\log(1/t)] = \lim_{t \to 0}\frac{\frac{d}{dt}\log[\log(1/t)]}{-1/t^2}=0$$

and therefore we obtain

(B.15) $$\int_0^x \log[\log(1/t)]dt = x\log[\log(1/x)]-li(x)$$

Since $\int_0^1 \log[\log(1/t)]dt = -\gamma$ we find that

(B.16) $$\lim_{x \to 1}[x\log[\log(1/x)]-li(x)]=-\gamma$$

Integration gives us

$$\int_0^x \left[\frac{1}{1-t}+\frac{1}{\log t}\right]dt = -\log(1-t)+li(t)\Big|_0^x$$

$$= -\log(1-x)+li(x)$$

We note from (B.2) that

$$\gamma = \int_0^1 \left[ \frac{1}{1-y} + \frac{1}{\log y} \right] dy$$

and therefore we have

(B.17) $\qquad \lim_{x \to 1}[-\log(1-x) + li(x)] = \gamma$

Therefore we see that

(B.18) $\qquad \lim_{x \to 1}[x \log[\log(1/x)] - \log(1-x)] = 0$

Multiplying (B.12) by $x$, upon taking the limit, we see that

$$\lim_{x \to 1}[x li(x) - x \log(-\log x)] = \gamma$$

Combining this with (B.18) we obtain

$$\lim_{x \to 1}[x li(x) - li(x) + \log(1-x) - x \log(-\log x)] = 0$$

and hence

(B.19) $\qquad \lim_{x \to 1}[(x-1) li(x)] = 0$

This may also be verified using L'Hôpital's rule

$$\lim_{x \to 1}[(x-1) li(x)] = \lim_{x \to 1} \frac{li(x)}{1/(x-1)} = -\lim_{x \to 1} \frac{(x-1)^2}{\log x} = -\lim_{x \to 1} 2x(x-1) = 0$$

The limit (B.19) in conjunction with (B.12) implies that

(B.20) $\qquad \lim_{x \to 1}[(x-1) \log(-\log x)] = 0$

□

We refer back to the Frullani integral (1.13)

$$\log u = \int_0^1 \left[ \frac{y^{u-1}}{\log y} - \frac{1}{\log y} \right] dy$$

where the third derivative results in

$$\frac{2}{u^3} = \int_0^1 y^{u-1} \log y^2 \, dy$$

We now let $u = n$ and make a summation to obtain the well-known integral

$$2\varsigma(3) = 2\sum_{n=1}^{\infty} \frac{1}{n^3} = \int_0^1 \frac{\log y^2}{1-y} dy$$

□

The following is extracted from Brede's dissertation [7]. In 1999 Coppo showed that there exists a polynomial $p_n(z)$ such that

$$x^n = \int_0^{\infty} p_n(x - \log z) e^{-z} dz$$

With $p_n(z) = \sum_{k=0}^{n} a_k z^k$ we have

$$\int_0^{\infty} p_n(x - \log z) e^{-z} dz = \int_0^{\infty} \sum_{k=0}^{n} a_k (x - \log z)^k e^{-z} dz$$

$$= \sum_{k=0}^{n} a_k \int_0^{\infty} \sum_{l=0}^{k} \binom{k}{l} x^l (-\log z)^{k-l} e^{-z} dz$$

$$= \sum_{k=0}^{n} a_k \sum_{l=0}^{k} \binom{k}{l} (-1)^{k-l} \Gamma^{(k-l)}(1) x^l$$

Hence we have

$$\int_0^{\infty} p_n(x - \log z) e^{-z} dz = \sum_{k=0}^{n} \left[ \sum_{l=0}^{n-k} (-1)^l \binom{k+l}{l} \Gamma^{(l)}(1) a_{k+l} \right] x^k$$

and we then select the coefficients $a_k$ such that

$$\sum_{l=0}^{n-k} (-1)^l \binom{k+l}{l} \Gamma^{(l)}(1) a_{k+l} = \begin{cases} 1 & \text{for } k = n \\ 0 & \text{for } k \text{ less than } n \end{cases}$$

resulting in

$$x^n = \int_0^\infty p_n(x - \log z) e^{-z} dz$$

For example we have

$$p_0(z) = 1$$

$$p_1(z) = z - \gamma$$

$$p_2(z) = z^2 - 2\gamma z + \gamma^2 - \varsigma(2)$$

Letting $x \to \log x$ we obtain

$$\log^n x = \int_0^\infty p_n(\log x - \log z) e^{-z} dz$$

and, with the substitution $z/x = -\log t$, we have

$$\frac{\log^n x}{x} = \int_0^1 p_n[-\log\log(1/t)] t^{x-1} dt$$

Referring to the well-known expression for the Stieltjes constants [19, p.4]

$$\gamma_n = \sum_{k=1}^\infty \left[ \frac{\log^n k}{k} - \frac{\log^n(k+1) - \log^n k}{n+1} \right] = \sum_{k=1}^\infty \left[ \frac{\log^n k}{k} - \int_k^{k+1} \frac{\log^n t}{t} dt \right]$$

$$= \lim_{N \to \infty} \left[ \sum_{k=1}^N \frac{\log^n k}{k} - \int_1^{N+1} \frac{\log^n t}{t} dt \right]$$

we then see that

$$\gamma_n = \lim_{N \to \infty} \left[ \sum_{k=1}^N \int_0^1 p_n[-\log\log(1/t)] t^{k-1} dt - \int_1^{N+1} \int_0^1 p_n[-\log\log(1/t)] t^{x-1} dt\, dx \right]$$

$$= \lim_{N \to \infty} \int_0^1 p_n[-\log\log(1/t)] \left[ \sum_{k=0}^{N-1} t^k - \int_1^{N+1} t^x dx \right] dt$$

$$= \lim_{N \to \infty} \int_0^1 p_n[-\log\log(1/t)] \left[ \frac{t^N - 1}{t - 1} - \frac{t^N - 1}{\log t} \right] dt$$

Hence, as shown more completely in Brede's dissertation, we obtain as $N \to \infty$

(B.21) $$\gamma_n = \int_0^1 P_n[-\log\log(1/t)] \left[\frac{1}{\log t} - \frac{1}{t-1}\right] dt$$

Since $\dfrac{1}{\log t} - \dfrac{1}{t-1} = \int_0^1 \dfrac{1-u^t}{1-u} du$ we may write this as a double integral

$$\gamma_n = \int_0^1 \int_0^1 P_n[-\log\log(1/t)] \frac{1-u^t}{1-u} du\, dt$$

$\square$

The following is extracted from Coppo's paper [14]. With the substitution $u = xt$ we see that

$$\int_0^\infty e^{-xt} \log^n t\, dt = \frac{1}{x} \int_0^\infty e^{-u} (\log u - \log x)^n du$$

$$= \sum_{k=0}^n \binom{n}{k} (-1)^k \frac{\log^k x}{x} \int_0^\infty e^{-u} \log^{n-k} u\, du$$

$$= \sum_{k=0}^n \binom{n}{k} (-1)^k \frac{\log^k x}{x} \Gamma^{(n-k)}(1)$$

Let us assume that

(B.22) $$\sum_{k=0}^n \binom{n}{k} (-1)^k \Gamma^{(n-k)}(1) \hat{a}_k(t) = \log^n t$$

where $\hat{a}_k(t)$ are to be determined. With $n = 0$ we have

$$\hat{a}_0(t) = 1$$

and $n = 1$ gives us

$$\Gamma^{(1)}(1) - \hat{a}_1(t) = \log t$$

We then multiply (B.22) by $e^{-xt}$ and integrate to obtain

$$\sum_{k=0}^{n}\binom{n}{k}(-1)^k \Gamma^{(n-k)}(1)\int_0^\infty e^{-xt}\hat{a}_k(t)\,dt = \int_0^\infty e^{-xt}\log^n t\,dt$$

We then have

(B.23) $$\sum_{k=0}^{n}\binom{n}{k}(-1)^k \Gamma^{(n-k)}(1)\left[\frac{\log^k x}{x} - \int_0^\infty e^{-xt}\hat{a}_k(t)\,dt\right] = 0$$

and Coppo deduces that

$$\frac{\log^k x}{x} = \int_0^\infty e^{-xt}\hat{a}_k(t)\,dt$$

It seems to me that this solution is only unique provided the coefficients appearing in (B.23) are linearly independent.

The associated Stieltjes functions $\hat{\gamma}_k(x)$ are defined by

$$\hat{\gamma}_k(x) = \int_0^\infty e^{-xt}\left(\frac{1}{1-e^{-t}} - \frac{1}{t}\right)\hat{a}_k(t)\,dt$$

and we have in particular

$$\hat{\gamma}_0(x) = \int_0^\infty e^{-xt}\left(\frac{1}{1-e^{-t}} - \frac{1}{t}\right)dt = \log x - \psi(x)$$

and with $x = 1$ we have

$$\hat{\gamma}_0(1) = \int_0^\infty e^{-t}\left(\frac{1}{1-e^{-t}} - \frac{1}{t}\right)dt = \gamma$$

It is easily seen that

(B.24) $$\sum_{k=0}^{n}\binom{n}{k}(-1)^k \Gamma^{(n-k)}(1)\hat{\gamma}_k(x) = \int_0^\infty e^{-xt}\left(\frac{1}{1-e^{-t}} - \frac{1}{t}\right)\log^n t\,dt$$

We have the well-known integral representation for the Hurwitz zeta function [26, p.92]

$$\varsigma(s,x) = \frac{1}{\Gamma(s)}\int_0^\infty \frac{t^{s-1}e^{-(x-1)t}}{e^t - 1}\,dt = \frac{1}{\Gamma(s)}\int_0^\infty \frac{t^{s-1}e^{-xt}}{1-e^{-t}}\,dt$$

and from the definition of the gamma function we see that

$$\frac{x^{1-s}}{s-1} = \frac{1}{\Gamma(s)} \int_0^\infty t^{s-2} e^{-xt} dt$$

Hence we have for $\mathrm{Re}(s) > 1$

$$\varsigma(s,x) - \frac{x^{1-s}}{s-1} = \frac{1}{\Gamma(s)} \int_0^\infty e^{-xt} \left( \frac{1}{1-e^{-t}} - \frac{1}{t} \right) t^{s-1} dt$$

which we write as

(B.25) $\quad \left( \varsigma(s,x) - \frac{x^{1-s}}{s-1} \right) \Gamma(s) = \int_0^\infty e^{-xt} \left( \frac{1}{1-e^{-t}} - \frac{1}{t} \right) t^{s-1} dt$

We have

$$\varsigma(s,x) - \frac{x^{1-s}}{s-1} = \varsigma(s,x) - \frac{1}{s-1} + \frac{1}{s-1} - \frac{x^{1-s}}{s-1}$$

and we denote $f(x)$ as

$$f(x) = \frac{1 - x^{1-s}}{s-1}$$

We can represent $f(x)$ by the following integral

$$f(x) = \frac{1 - x^{1-s}}{s-1} = \int_1^x t^{-s} dt$$

so that

$$f^{(r)}(x) = (-1)^r \int_1^x t^{-s} \log^r t \, dt$$

and thus

$$f^{(r)}(1) = (-1)^r \int_1^x \frac{\log^r t}{t} dt = (-1)^r \frac{\log^{r+1} x}{r+1}$$

Differentiating (B.25) $n$ times gives us

$$\frac{d^n}{ds^n} \left( \varsigma(s,x) - \frac{x^{1-s}}{s-1} \right) \Gamma(s) = \int_0^\infty e^{-xt} \left( \frac{1}{1-e^{-t}} - \frac{1}{t} \right) t^{s-1} \log^n t \, dt$$

The Stieltjes constants $\gamma_n(u)$ are the coefficients in the Laurent expansion of the Hurwitz zeta function $\varsigma(s,u)$ about $s=1$

$$\varsigma(s,x) = \frac{1}{s-1} + \sum_{k=0}^{\infty} \frac{(-1)^k}{k!} \gamma_k(x)(s-1)^k$$

and we therefore have

$$\frac{d^r}{ds^r}\left[\varsigma(s,x) - \frac{1}{s-1}\right] = \sum_{k=0}^{\infty} \frac{(-1)^k}{k!} \gamma_k(x) k(k-1)\cdots(k-r+1)(s-1)^{k-r}$$

and thus

$$\lim_{s\to 1} \frac{d^r}{ds^r}\left[\varsigma(s,x) - \frac{1}{s-1}\right] = (-1)^r \gamma_r(x)$$

$$\lim_{s\to 1} \frac{d^r}{ds^r}\left[\varsigma(s,x) - \frac{x^{1-s}}{s-1}\right] = (-1)^r \gamma_r(x) + (-1)^r \frac{\log^{r+1} x}{r+1}$$

Applying the Leibniz rule we have

$$\lim_{s\to 1} \frac{d^n}{ds^n}\left(\varsigma(s,x) - \frac{x^{1-s}}{s-1}\right)\Gamma(s) = \sum_{k=0}^{n} \binom{n}{k}(-1)^k \Gamma^{(n-k)}(1)\left[\gamma_k(x) + \frac{\log^{k+1} x}{k+1}\right]$$

and hence using (B.24) we see that

(B.26) $$\sum_{k=0}^{n}\binom{n}{k}(-1)^k \Gamma^{(n-k)}(1)\left[\gamma_k(x) + \frac{\log^{k+1} x}{k+1}\right] = \int_0^{\infty} e^{-xt}\left(\frac{1}{1-e^{-t}} - \frac{1}{t}\right)\log^n t\, dt$$

We then see from (B.24) that (subject to the previous comment regarding the coefficients being linearly independent)

$$\hat{\gamma}_k(x) = \gamma_k(x) + \frac{\log^{k+1} x}{k+1}$$

and accordingly that

$$\hat{\gamma}_k(1) = \gamma_k(1) = \gamma_k$$

With $x=1$ we have

(B.27) $$\sum_{k=0}^{n}\binom{n}{k}(-1)^k \Gamma^{(n-k)}(1)\gamma_k = \int_0^\infty e^{-t}\left(\frac{1}{1-e^{-t}}-\frac{1}{t}\right)\log^n t\, dt$$

We also have

$$\varsigma(s,x) - \frac{x^{1-s}}{s-1} = \sum_{k=0}^{\infty}\frac{(-1)^k}{k!}\hat{\gamma}_k(x)(s-1)^k$$

With the substitution $t = \log(1/u)$ we obtain

(B.28) $$\sum_{k=0}^{n}\binom{n}{k}(-1)^k \Gamma^{(n-k)}(1)\hat{\gamma}_k(x) = \int_0^1 u^{x-1}\left(\frac{1}{1-u}+\frac{1}{\log u}\right)\log^n \log\left(\frac{1}{u}\right)du$$

and

(B.29) $$\sum_{k=0}^{n}\binom{n}{k}(-1)^k \Gamma^{(n-k)}(1)\gamma_k = \int_0^1 \left(\frac{1}{1-u}+\frac{1}{\log u}\right)\log^n \log\left(\frac{1}{u}\right)du$$

and in particular we have

$$A_0(x) = \hat{\gamma}_0(x) = \log x - \psi(x)$$

$$A_1(x) = \Gamma^{(1)}(1)\hat{\gamma}_0(x) - \hat{\gamma}_1(x)$$

This gives us

$$\hat{\gamma}_1(x) = -\gamma[\log x - \psi(x)] - \int_0^1 u^{x-1}\left(\frac{1}{1-u}+\frac{1}{\log u}\right)\log\log\left(\frac{1}{u}\right)du$$

$$\gamma_1 = -\gamma^2 - \int_0^1 \left(\frac{1}{1-u}+\frac{1}{\log u}\right)\log\log\left(\frac{1}{u}\right)du$$

For $u \in [0,1]$ we have

$$\log\left(\frac{1}{u}\right) = \left|\log\left(\frac{1}{u}\right)\right| = |-\log u| = |\log u|$$

and hence we obtain

(B.30) $$\gamma_1 = -\gamma^2 - \int_0^1 \left( \frac{1}{1-u} + \frac{1}{\log u} \right) \log|\log u| \, du$$

which unfortunately differs from (2.5).

## Appendix C

## Derivatives of the gamma function

The gamma function is defined as

$$\Gamma(x) = \int_0^\infty y^{x-1} e^{-y} dy$$

and, with the substitution $y = -\log t$, this becomes

$$\Gamma(x) = \int_0^1 \left( \log(1/t) \right)^{x-1} dt$$

We therefore have

$$\Gamma'(x) = \int_0^1 \left( \log(1/t) \right)^{x-1} \log\log(1/t) \, dt$$

and thus

(C.1) $$\Gamma'(1) = \int_0^1 \log\log(1/t) \, dt = -\gamma$$

More generally we see that

$$\Gamma^{(n)}(x) = \int_0^1 \log^{x-1}(1/t) \left( \log\log(1/t) \right)^n dt$$

and as we shall see below

$$\Gamma^{(n)}(1) = \int_0^1 \left( \log\log(1/t) \right)^n dt = (-1)^n c_n$$

where $c_n$ are positive constants.

We have

$$\Gamma^{(n)}(x) = \int_0^\infty e^{-t} t^{x-1} \log^n t \, dt$$

and thus

$$\Gamma^{(n)}(1) = \int_0^\infty e^{-t} \log^n t \, dt$$

We see that

$$\int_0^\infty e^{-t} \log^n t \, dt = \int_0^1 e^{-t} \log^n t \, dt + \int_1^\infty e^{-t} \log^n t \, dt$$

and with the substitution $t = 1/y$ we have

$$\int_1^\infty e^{-t} \log^n t \, dt = (-1)^n \int_0^1 e^{-1/y} y^{-2} \log^n y \, dy$$

We therefore obtain

$$\Gamma^{(n)}(1) = \int_0^1 \left[ t^2 + (-1)^n e^{t-1/t} \right] e^{-t} t^{-2} \log^n t \, dt$$

and it is clear that the integrand is positive in the interval $[0,1]$ in the case where $n$ is an even integer. We now consider the case where $n$ is an odd integer and we want to prove that for $t \in [0,1]$

$$t^2 \geq e^{t-1/t}$$

and this is easily seen by considering the logarithmic equivalent

$$2 \log t \geq t - 1/t$$

This inequality is easily proved by considering

$$f(t) = 2t \log t - t^2 + 1$$

where $f(0) = 1$ and $f'(t) = 2[\log t + 1 - t]$ and it is known that $\log t > 1 - t$ for $t > 0$.

Therefore we see that $\Gamma^{(n)}(1)$ has the same sign as $(-1)^n$. This result was reported as an exercise in Apostol's book [4, p.303] and we used it in (3.13) above.

**REFERENCES**

[1] V.S.Adamchik, A Class of Logarithmic Integrals. Proceedings of the 1997
    International Symposium on Symbolic and Algebraic Computation.
    ACM, Academic Press, 1-8, 2001.
    http://www-2.cs.cmu.edu/~adamchik/articles/logs.htm

[2] O.R. Ainsworth and L.W. Howell, The generalized Euler-Mascheroni constants.
    NASA Centre for AeroSpace Information (CASI)
    NASA-TP-2264 ;NAS 1.60.2264, 1984. View PDF File

[3] J. Anglesio, A fairly general family of integrals.
    Amer. Math. Monthly, 104, 665-666, 1997.

[4] T.M. Apostol, Mathematical Analysis, Second Ed., Addison-Wesley Publishing
    Company, Menlo Park (California), London and Don Mills (Ontario), 1974.

[5] E.T. Bell, Exponential polynomials. Ann. of Math., 35 (1934), 258-277.

[6] G. Boros and V.H. Moll, Irresistible Integrals: Symbolics, Analysis and
    Experiments in the Evaluation of Integrals. Cambridge University Press, 2004.

[7] M. Brede, Eine reihenentwicklung der vervollständigten und ergänzten
    Riemannschen zetafunktion und verwandtes
    http://www.mathematik.uni-kassel.de/~koepf/Diplome/brede.pdf

[8] T.J.I'a Bromwich, Introduction to the theory of infinite series. 2$^{nd}$ edition
    Macmillan & Co Ltd, 1965.

[9] B. Candelpergher, M.A. Coppo and E. Delabaere, La Sommation de
    Ramanujan. L'Enseignement Mathématique, 43, 93-132, 1997.
    http://math.univ-angers.fr/~delab/Ramanujan.ps

[9a] M.W. Coffey, New results concerning power series expansions of the Riemann xi
    function and the Li/Keiper constants. (preprint) 2007.

[9b] M.W. Coffey, Relations and positivity results for the derivatives of the Riemann
    $\xi$ function. J. Comput. Appl. Math., 166, 525-534 (2004)

[10] C.B. Collins, The role of Bell polynomials in integration.
    J. Comput. Appl. Math. 131 (2001) 195-211.

[11] L. Comtet, Advanced Combinatorics, Reidel, Dordrecht, 1974.

[12] D.F. Connon, Some series and integrals involving the Riemann zeta function,
    binomial coefficients and the harmonic numbers. Volume II(a), 2007.

       arXiv:0710.4023 [pdf]

[13] D.F. Connon, Some series and integrals involving the Riemann zeta function,
     binomial coefficients and the harmonic numbers. Volume II(b), 2007.
       arXiv:0710.4024 [pdf]

[14] M.A. Coppo, Nouvelles expressions des constantes de Stieltjes.
     Expositiones Mathematicae 17, No. 4, 349-358 (1999).

[15] B. Fuglede, A sharpening of Wielandt's characterization of the gamma function.
     Amer. Math. Monthly, 115, 845-850, 2008.

[16] I.S. Gradshteyn and I.M. Ryzhik, Tables of Integrals, Series and Products.
     Sixth Ed., Academic Press, 2000.
     Errata for Sixth Edition http://www.mathtable.com/errata/gr6_errata.pdf

[17] J. Guillera and J. Sondow, Double integrals and infinite products for some
     classical constants via analytic continuations of Lerch's transcendent. 2005.
     math.NT/0506319 [abs, ps, pdf, other]

[18] H. Hasse, Ein Summierungsverfahren für Die Riemannsche $\varsigma$ - Reithe.
     Math.Z.32, 458-464, 1930.
       http://dz-srv1.sub.uni-goettingen.de/sub/digbib/loader?ht=VIEW&did=D23956&p=462

[19] A. Ivić, The Riemann Zeta- Function: Theory and Applications.
     Dover Publications Inc, 2003.

[20] K.S. Kölbig, The complete Bell polynomials for certain arguments in terms of
     Stirling numbers of the first kind.
     J. Comput. Appl. Math. 51 (1994) 113-116.
     Also available electronically at:
     A relation between the Bell polynomials at certain arguments and a Pochhammer
     symbol. CERN/Computing and Networks Division, CN/93/2, 1993.
     http://doc.cern.ch/archive/electronic/other/preprints//CM-P/CM-P00065731.pdf

[21] N. Nielsen,Theorie des Integrallogarithmus und verwanter tranzendenten.1906.
       http://www.math.uni-bielefeld.de/~rehmann/DML/dml_links_author_H.html

[22] N. Nielsen, Die Gammafunktion. Chelsea Publishing Company, Bronx and
      New York, 1965.

[23] S.K. Lakshamana Rao, On the Sequence for Euler's Constant.
     Amer. Math. Monthly, 63, 572-573, 1956.

[24] T. Rivoal, Polynômes de type Legendre et approximations de la constant
      d'Euler. Note (2005).   DVI,   PS,   PDF

http://www-fourier.ujf-grenoble.fr/

[25] G.K. Srinivasan, The gamma function: An eclectic tour.
Amer. Math. Monthly, 114, 297-315, 2007.

[26] H.M. Srivastava and J. Choi, Series Associated with the Zeta and Related Functions. Kluwer Academic Publishers, Dordrecht, the Netherlands, 2001.

[27] E.C. Titchmarsh, The Theory of Functions.
2$^{nd}$ Ed., Oxford University Press, 1932.

Donal F. Connon
Elmhurst
Dundle Road
Matfield
Kent TN12 7HD
England

dconnon@btopenworld.com